\documentclass[12pt]{article}
\usepackage{latexsym}
\usepackage{amssymb}
\def\theequation{\thesection.\arabic{equation}}

\font\gorditas = msbm8
\def\bbb#1{\hbox {{\gordas #1}}}

\def\enita{\hbox{\gorditas N}}
\def\zetita{\hbox{\gorditas Z}}

\font\gordas = msbm10 at 12pt
\def\bbb#1{\hbox {{\gordas #1}}}
\def\erre{{\bbb R}}
\def\zet{{\bbb Z}}

\def\UNO{1\mkern-7mu1}

\newtheorem{theorem}{Theorem}[section]
\newtheorem{lemma}[theorem]{Lemma}

\newtheorem{proposition}[theorem]{Proposition}
\newtheorem{ex}[theorem]{Example}

\newtheorem{remark}[theorem]{Remark}

\begin{document}
\noindent

\begin{center}
{\bf\large Number variance for hierarchical random walks\\
and related fluctuations$^*$}
 \vglue.5cm

\begin{tabular}{ccc}
Tomasz Bojdecki$^1$& Luis G. Gorostiza$^\dagger$&Anna  Talarczyk$^1$\\
tobojd@mimuw.edu.pl&lgorosti@math.cinvestav.mx&annatal@mimuw.edu.pl
\end{tabular}
\end{center}
\vglue.3cm
\begin{abstract} 
We study an infinite system of independent symmetric random walks on a hierarchical group, in particular, the $c$-random walks. Such walks are used, e.g., in population genetics. The number variance problem consists in investigating if the variance of the number of ``particles'' $N_n(L)$ lying in the ball of radius $L$ at a given time $n$ remains bounded, or even better,  converges to a finite limit, as $L\to \infty$. We give a necessary and sufficient condition and discuss its relationship to transience/recurrence property of the walk. Next we consider normalized fluctuations of $N_n(L)$ around the mean as $n\to \infty$ and $L$ is increased in an appropriate way. We prove convergence of finite dimensional distributions to a Gaussian process whose properties are discussed. As the $c$-random walks mimic symmetric stable processes on $\erre$, we compare our results to those obtained by Hambly and Jones (2007,2009), where the number variance problem for an infinite system of symmetric stable processes on $\erre$ was studied.
 Since the
hierarchical group is an ultrametric space, corresponding results for
symmetric stable processes and hierarchical random
walks may be analogous or quite different, as has been observed in
other contexts.
An example of a difference in the present context is that for the stable processes a fluctuation limit process is a centered Gaussian process which is not Markovian and has 
long
range dependent stationary increments, but the counterpart for hierarchical
random walks is Markovian, and in a special case it has
independent increments.
\end{abstract}
\vglue1cm
\noindent
{\bf AMS 2000 subject classifications:} Primary 60G50, Secondary 60F05.\\
{\bf Key words:} hierarchical random walk, hierarchical group, ultrametric, number variance, fluctuation, limit theorem.

\footnote{\kern-.66cm $*$ Supported in part by CONACyT grant 98998 (Mexico) and MNiSzW grant N N201 397537 (Poland).\\
$^1$ Institute of Mathematics, University of Warsaw, ul. Banacha 2, 02-097 Warszawa, Poland\\
$^\dagger$ Centro de Investigaci\'on y de Estudios Avanzados, A.P. 14-740, Mexico 07000 D.F., Mexico.}

\section{Introduction}
\setcounter{section}{1}
\setcounter{equation}{0}

Random walks on hierarchical groups, called {\it hierarchical random walks}, have arisen in several applications. The hierarchical random walks we will 
consider in this paper are of
 the type introduced by Spitzer \cite{S} (p.93) in a special case, and more generally by Sawyer and Felsenstein \cite{SF} in  the context of population genetics.
We will focus on  ``number variance'' properties and related fluctuations   for these random walks. This  question has been investigated by Hambly and Jones \cite{HJ1,HJ2} for $\alpha$-stable processes on $\erre$.

For an integer $M\geq 2$, the {\it hierarchical group of order} $M$ is defined by
$$\Omega_M=\{x=(x_1,x_2,\ldots):x_i\in\{0,1,\ldots,M-1\},
\Sigma_ix_i<\infty\},$$
with addition componentwise mod $M$. It is a countable Abelian group which is also described as the direct sum of a countable number of copies of the cyclic group of order $M$. The {\it hierarchical distance} $|\cdot|$ on $\Omega_N$ is defined by
$$|x-y|=\left\{
\begin{array}{lll}
0&{\rm if}&x=y,\\
\max\{i:x_i\neq y_i\}&{\rm if}&x\neq y.
\end{array}\right.
$$
It is translation-invariant, and satisfies the strong (non-Archimedean) triangle inequality
$$|x-y|\leq\max\{|x-z|,|z-y|\}\quad\hbox{\rm for any}\quad x,y,z.$$
This means that $(\Omega_M,|\cdot|)$ is an ultrametric space. In such a space two balls are either disjoint or one is contained in the other. This property is the cause of significant differences from mathematical models based on Euclidean spaces $\erre^d$ or Euclidean lattices $\zet^d$, and it is interesting to study analogies and differences between ``corresponding'' models defined on $\Omega_M$ and on $\erre^d$ or $\zet^d$, specially because hierarchical models may reveal behaviors that are not observed in Euclidean  models. A picture of $(\Omega_M,|\cdot|)$ is the set of leaves at the top of an infinite tree where each inner node at each level $j\geq 1$ has one predecessor at level $j+1$ and $M$ successors at level $j-1$. The {\it distance} between two individuals (leaves) at level $j=0$ is the depth in the tree to their most recent common ancestor, and it measures the degree of relatedness between the two individuals.

A random walk $\xi=(\xi_n)_{n\geq 0}$
 on $\Omega_M$ starting at $0$ is defined by $\xi_0=0,\xi_n=\rho_1+\cdots+\rho_n,n\geq 1$, where $\rho_1,\rho_2,\ldots$ are independent copies of $\rho$, which is a random element of $\Omega_M$ with distribution of the form
\begin{equation}
\label{eq:1.1}
P(\rho=y)=\frac{r_{|y|}}{M^{|y|-1}(M-1)},\quad y\neq 0,\quad P(\rho=0)=0,
\end{equation}
where $(r_j)_{j=1,2,\ldots}$ is a probability law on $\{1,2,\ldots\}$. That is, the jumps of $\xi_n$ are taken by first choosing distance $j$ with probability 
$r_j$, and then choosing a point with uniform probability among those at distance $j$ from the previous position of the walk (note that $M^{j-1}(M-1)$ is the number of points at distance $j$ from a given point of $\Omega_M$). These random walks are the most general symmetric random walks on $\Omega_M$. We assume that $(r_j)$ is not restricted to a bounded set, so that the walk can reach arbitrarily large distances (by ultrametricity, it is not possible to go far with small steps). We refer to a hierarchical random walk determined by $(r_j)$ as 
$r_j$-{\it rw}. A particular $r_j$-rw is the 
{$c$-{\it rw}}, where $c$ is a constant such that $0<c<M$, and
\begin{equation}
\label{eq:1.2}
r_j=\left(1-\frac{c}{M}\right)\left(\frac{c}{M}\right)^{j-1}, \quad j=1,2,\ldots,
\end{equation}
The $c$-rw mimics the behavior of the standard $\alpha$-stable process on $\erre^d$  in the sense that both have analogous recurrence/transience behaviors. The analogies are given in terms of their degrees  $\gamma$. The Appendix contains  information on $\gamma$. In particular, $\gamma<0$ corresponds to recurrence and $\gamma>0$ to transience. The $\alpha$-stable process on $\erre^d$ has degree
\begin{equation}
\label{eq:1.3}
\gamma=\frac{d}{\alpha}-1,
\end{equation}
and the $c$-rw on $\Omega_M$ has degree
\begin{equation}
\label{eq:1.4}
\gamma=\frac{\log c}{\log(M/c)}.
\end{equation}
By relating the degrees (\ref{eq:1.3}) and (\ref{eq:1.4}) we will compare our results with those of \cite{HJ1,HJ2}.
Other  $r_j$-rw's we will consider
are defined   in the Appendix   and used for examples. 

Hierarchical structures are found in the physical, biological and social sciences due to the multiscale organization of many natural objects (see e.g. 
\cite{BF,RTV}). Some of the main applications of hierarchical systems are found in statistical physics; a partial list of references is contained in \cite{DGW3}. We have already mentioned \cite{SF} in population genetics. Some other references where stochastic models based on hierarchical structures have been 
studied  are 
\cite{DGr,DGHSS,FG,K} 
(interacting diffusions), \cite{DGW1,DGW4} (branching systems), \cite{AS} (contact processes), \cite{DG1,DG2,KMT} (percolation), \cite{K1,K2} (search algorithms).

The ``number variance'' problem studied  \cite{HJ1,HJ2} for $\alpha$-stable processes on $\erre$ has a motivation from physics related to so-called determinantal processes (see the Introduction of \cite{HJ1}). 
The model is  as follows in  simplified form. An independent copy of 
the standard $\alpha$-stable process, representing the motion of a particle, starts at time $0$ from each point $u_j=j-\varepsilon, j\in\zet$, where $\varepsilon$ is a uniformly distributed random variable on [0, 1] (the $\varepsilon$-displacement provides  spatial homogeneity).  Let $N_t[0,L]$ denote the number of particles that at time $t>0$ lie in the interval $[0,L]$, and consider its variance,  
Var$N_t[0,L]$. The ``number variance'' problem refers to the behavior of 
Var$N_t[0,L]$ for fixed $t$ as $L\to\infty$. It is shown that for $\alpha<1$ (transient case) the variance tends to $\infty$, and for $1<\alpha\leq 2$ (recurrent case) it has a finite limit up to an additive fluctuation bounded by 1. The latter behavior is called {\it saturation}. For $\alpha=1$ the process is recurrent and the variance tends to infinity. The rescaled fluctuation process defined by
\begin{equation}
\label{eq:1.5}
Z_t(s)=\frac{N_t[0,st^{1/\alpha}]-EN_t[0,st^{1/\alpha}]}{t^{1/2\alpha}},\quad
 s\geq 0,
\end{equation}
is discussed, and a continuous interpolation of this process is shown to converge weakly on a space of continuous functions,  as $t\to\infty$, to a centered Gaussian process with stationary increments, which is not Markovian and has long range dependence. Hence this process has a resemblance to fractional Brownian motion.

The purpose of the present paper is to prove counterparts to the results of \cite{HJ1,HJ2} for a system of independent hierarchical random walks starting from each point of $\Omega_M$. We will do this for a  general class of hierarchical random walks, and we will also consider a model where at the initial time there is a Poisson number of particles at each site. As with other systems of hierarchical random walks, one of our objectives is to exhibit analogies and differences with the results  for the Euclidean  model.

Now we give a summary of our results.

We consider a system of independent $r_j$-rw's (particles) starting from each point of $\Omega_M$, and denote by $N_n (L)$ the number of particles that at step $n$ lie in the ball $\{x\in \Omega_M: |x|\leq L\}$, 
$L\in \zet_+$ (in \cite{HJ1,HJ2} the interval $[0,L]$ was considered, but essentially nothing changes if it is replaced by $[-L, L]$). For general $r_j$ we give a necessary and sufficient condition for boundedness of Var$N_n(L)$ and for existence of its limit as $L\rightarrow \infty$. It turns out that
$$
\lim_{L\rightarrow \infty} {\rm Var}N_n (L)=2n \frac{1}{M-1} \lim_{L\rightarrow \infty} M^Lr_L
$$
(Theorem 2.2). In particular, for the $c$-rw this limit exists, and it is finite if and only if $c\leq 1$, i.e., if the walk is recurrent. For the non-critical cases, $c<1$ and $c>1$, this result corresponds to the properties of the number variance for $\alpha$-stable processes mentioned above, but the correspondence breaks down for the critical case, $c=1$: the process is recurrent and the  variance has a finite limit. For general  $r_j$, boundedness of Var$N_n (L)$ and recurrence are not equivalent (see Remark 2.4(c)). For the initial Poisson system the situation is simpler, 
Var$N_n (L)$ always tends to infinity. It can be shown that the same thing happens for the $\alpha$-stable process.

Next we investigate the fluctuations of $N_n (L)$ around the mean as $n\rightarrow \infty$, when $L=L(n)$ is increased in an appropriate way. Analogously as in (\ref{eq:1.5}), where a ``time'' parameter $s$ was introduced, we introduce a ``time'' parameter $t$, taking $L(n)+R(t)$, where $R(t)$ is a suitable non-decreasing function of $t$, and we consider a corresponding fluctuation process. We assume that 
$\lim_{j\rightarrow \infty}r_{j+1}/r_j=a>0$. It
turns out that the cases $a<1$ and $a=1$ are significantly different (see Theorem 2.8). For $a<1$, the most appropriate $L(n)$ has logarithmic growth and the norming is of the order $n^{\theta /2}$, where $\theta =\log M/\log \frac{1}{a}$. In general the fluctuation processes may not converge, but convergent subsequences can be chosen (in the sense of finite dimensional distributions), and a family of processes indexed by a parameter $\kappa \in [1, \frac{1}{a}]$ is obtained. These processes have stepwise trajectories, they are Gaussian with covariance of the form $f(s\vee t)g_\kappa (s\wedge t)$, in particular they are Markov. If $a=1$, the fluctuation process converges to a Gaussian process with independent increments. In contrast, for the $\alpha$-stable process on $\mathbb{R}$ the fluctuation limit is a non-Markov process with stationary increments \cite{HJ1,HJ2}.

We study some further properties of our limit processes and discuss a probabilistic interpretation of the parameter $\theta$, which is closely related to the degree $\gamma$.

Finally, for the Poisson model we show that with general $r_j$ the fluctuation limit exists and it is, up to a constant, the same process as that  obtained in the previous case for $a=1$.

Due to the ultrametric structure of $\Omega_M$, the calculations involved in the proofs are quite different from those in \cite{HJ1, HJ2}.
\section{Results}
\setcounter{section}{2}
\setcounter{equation}{0}

Fix an integer $M\geq 2$ and let $\Omega_M$ be the corresponding hierarchical group as defined in Introduction. We consider a system of independent $r_j$-rw's on $\Omega_M$ 
 starting from each point $u\in\Omega_M$.
Recall that $N_n(L)$ is the number of particles lying in the ball 
$B_L=\{x\in\Omega_M:|x|\leq L\}$ at step $n$, i.e.,
\vglue.3cm
\noindent
\begin{equation}
\label{eq:2.1a}
N_n(L)=\sum\limits_{u\in\Omega_M}\UNO_{B_L}(u+\xi_n^u),
\end{equation}
where $\{\xi^u\}_{u\in\Omega_M}$ are independent copies of $r_j$-rw's starting at 0.
\begin{proposition}
\label{p:2.1}

\begin{eqnarray}
\label{eq:2.1}
 &\kern -3.5cm (a)& EN_n(L)=M^L,\,\, n,\,\, L=0,1,2,\ldots,\\
\label{eq:2.2}
&\kern -3.5cm(b)& {\rm Var} N_n(L)\sim 2n M^LP(|\rho|>L)\,\,\hbox{as}\,\, L\to\infty,\\
\label{eq:2.3}
&\kern -3.5cm(c)& P(|\rho|\leq L)^{2n-1}nP(|\rho|>L)M^L\leq {\rm Var}N_n(L)\leq 2M^L,\quad 
 n,\,\, L=0,1,2,\ldots,
\end{eqnarray}
{\rm($\sim$ means that the quotient of both sides tends to 1).}
\end{proposition}

The number variance problem is solved in the following theorem.

\begin{theorem}
\label{t:2.2}
For any $n=1,2,\ldots, \lim\sup_{L\to\infty}{\rm Var} N_n(L)<\infty$ if and only if \\
$\lim\sup_{L\to\infty}M^Lr_L<\infty$. Moreover, 
$\lim_{L\to\infty}{\rm Var}N_n(L)$ exists (finite or not) if and only if $\lim_{L\to\infty}M^Lr_L$ exists, and in this case
\begin{equation}
\label{eq:2.4}
\lim_{L\to\infty}{\rm Var}N_n(L)=2n\frac{1}{M-1}\lim_{L\to\infty}
M^L{r_L}.
\end{equation}
\end{theorem}

This theorem applied to the examples mentioned
 in the Appendix gives the following results.

\begin{ex}
\label{2.3}

(a) $c$-rw:
\begin{equation}
\label{eq:2.5}
\lim_{L\to\infty}{\rm Var}N_n(L)=\left\{
\begin{array}{lll}
0&{\it if}&c<1,\nonumber\\
2n&{\it if}& c=1,\\
\infty&{\it if}&c>1.\nonumber
\end{array}\right.
\end{equation}
(b) $j^\beta$-rw:
$$\lim_{L\to\infty}{\rm Var}N_n(L)=\left\{
\begin{array}{lll}
\infty&{\it if}&\beta>0,\\
2n&{\it if}& \beta=0.
\end{array}\right.
$$
(for $\beta=0$ it is obviously a $c$-rw with $c=1$)

\noindent
(c) $r_j=Dj^{-\beta},\beta>1$:
$$\lim_{L\to\infty}{\rm Var}N_n(L)=\infty.$$
\end{ex}

\begin{remark}
\label{r:2.4} {\rm (a) Let us compare the result (\ref{eq:2.5}) for the $c$-rw with the solution of the number variance problem obtained in [13,14] for the $\alpha$-stable process. With $c\leq 1$ (recurrent case) we obtain ``true'' limits, whereas for the $\alpha$-stable process with $\alpha>1$ the variances are bounded but do not converge due to an oscillating term (see [14]). On the other hand, with  $c<1$ the limit is trivial (zero), whereas for the $\alpha$-stable process with $\alpha>1$ the variances are bounded away from $0$. The only non-trivial limit for the $c$-rw is obtained in the critical case $c=1$, and it has no counterpart for the $\alpha$-stable process, since for $\alpha=1$ the variances tend to infinity.

\noindent
(b) For the $c$-rw, finiteness of the variance limit is equivalent to recurrence, whereas for the $\alpha$-stable process this equivalence breaks down in the critical case $\alpha=1$.

\noindent
(c) Example 2.3(b) shows that in general finiteness of the variance limit is not equivalent to recurrence (see the Appendix).}
\end{remark}

We now give a result for the Poisson case.

\begin{proposition}
\label{p:2.5}
Assume that initially at each site there is a Poisson number of particles, and these numbers are i.i.d. Then for any system of $r_j$-rw's, 
$\lim_{L\to\infty}{\rm Var} N_n(L) = \infty$.
\end{proposition}

\begin{remark}
\label{r:2.6} {\rm (a) It can be shown that an analogous result holds for $\alpha$-stable processes.

\noindent
(b) We will come back to the Poisson system later on.  Poisson systems seem to be the most natural as random initial configurations. However, it can be shown that for each initial configuration determined by i.i.d. random variables $\{\nu_x\}_{x\in\Omega_M}(\nu_x$ particles at site $x$) which are truly random, i.e., with
Var$\;\nu_x>0$, we have $\lim_{L\to\infty}{\rm Var}N_n(L)=\infty$.}
\end{remark}

So far we have considered $N_n(L)$ for $n$ fixed. Now we will vary both $n$ and $L$, more precisely, we want to investigate the normalized fluctuations of $N_n(L)$ as $n\to\infty$, simultaneously increasing $L$ in an appropriate way.

We make the following assumption on the random walk: 
\begin{equation}
\label{eq:2.6}
r_{j+1}\leq r_j\,\,\hbox{\rm for}\,\, j\geq\,\,{\rm some}\,\, j_0,\,\,\hbox{\rm and}\,\, 
\lim_{j\to\infty}\frac{r_{j+1}}{r_j}=a>0.
\end{equation}
This assumption is satisfied for all our examples (see Example 2.3 and the 
Appendix).

By Proposition 2.1(c), it is natural to take
\begin{equation}
\label{eq:2.7}
L(n)=\sup\left\{L\in\zet_+:h(L)\geq \frac{1}{n}\right\},
\end{equation}
where 
\begin{equation}
\label{eq:2.8}
h(L)=P(|\rho|>L)=\sum^\infty_{j=L+1}r_j,
\end{equation}
so that, by (\ref{eq:2.3}), Var$N_n(L(n))$ has the same rate of increase as 
$M^{L(n)}$.  Hence 
$\sqrt{M^{L(n)}}$ is the natural normalization for the fluctuations of $N_n(L(n))$.

$L(n)$ has the following properties.

\begin{lemma}
\label{l:2.7} (a) if $a<1$, then
\begin{equation}
\label{eq:2.9}
\lim_{n\to\infty}\frac{L(n)}{\log_a\frac{1}{n}}=1.
\end{equation}
(b)
\begin{equation}
\label{eq:2.10}
1\leq n h(L(n))\quad\,\,\hbox{ and}\quad\,\,\lim\sup_{\kern-.8cm n\to\infty} n h(L(n))\leq\frac{1}{a}.
\end{equation}
\end{lemma}

As in [13,14], we want to investigate a fluctuation process introducing a new ``time'' parameter. In [13,14] this was done with a multiplicative parameter, but in our case an additive parameter is more adequate. This is caused by the hierarchical structure of the state space, which for $a<1$ implies the logarithmic growth of the radius of the balls (Lemma 2.7(a)). Consider any non-decreasing function $R:\erre_+\to\zet$ such that $\lim_{t\to\infty}R(t)=\infty$. For $n=1,2,\ldots$ we define the fluctuation process as
\begin{equation}
\label{eq:2.11}
X_n(t)=\frac{N_n((L(n)+R(t))^+)-E N_n((L(n)+R(t))^+)}{\sqrt{M^{L(n)}}},
\end{equation}
(see (\ref{eq:2.1})).

In what follows $\Rightarrow_f$ denotes weak convergence of finite-dimensional distributions, and $\lfloor x\rfloor$ is the integer part of $x\in\erre$. 

For $a$ given by
(\ref{eq:2.6}), we denote 
\begin{equation}
\label{eq:2.12}
b=\frac{M-a}{M-1}.
\end{equation}

\begin{theorem}
\label{t:2.8}
Assume (\ref{eq:2.6}).

\noindent
(a) Suppose $a<1$. Let $(n_i)_i$ be any subsequence such that
\begin{equation}
\label{eq:2.13}
\lim_{i\to\infty}n_i h(L(n_i))=\kappa
\end{equation}
for some $\kappa\in [1,\frac{1}{a}]$. Then $X_{n_i}\Rightarrow_f X^{(\kappa,R)}$, where $X^{(\kappa,R)}$ is a centered Gaussian process with covariance
\begin{equation}
\label{eq:2.14}
EX^{(\kappa,R)}(s)X^{(\kappa,R)}(t)=M^{R(s\wedge t)}g_\kappa(s\vee t),
\end{equation}
where
\begin{eqnarray}
g_\kappa(t)&=&1-\left(\frac{M-1}{M}\right)^2\left(\sum^\infty_{j=0}
\frac{e^{-\kappa ba^{R(t)}a^j}}{M^j}\right)^2\nonumber\\
\label{eq:2.15}
&-&\frac{(M-1)^3}{M^4}\sum^\infty_{j=0}\frac{1}{M^j}\left(\sum^\infty_{k=0}
\frac{e^{-\kappa ba^{R(t)}a^{j+k+1}}-
e^{-\kappa ba^{R(t)}a^j}}{M^k}\right)^2.
\end{eqnarray}
(b) If $a=1$, then $\lim_{n\to\infty} n h(L(n))=1$ and $X_n\Rightarrow_f X^{(R)}$, where $X^{(R)}$ is a centered Gaussian process with covariance
\begin{equation}
\label{eq:2.16}
EX^{(R)}(s)X^{(R)}(t)=(1-e^{-2})M^{R(s\wedge t)}.
\end{equation}
\end{theorem}

\begin{remark}
\label{r:2.9}

{\rm (a) Existence of a subsequence $(n_i)_i$ in part (a) follows immediately from Lemma 2.7(b). For example, for $c$-rw we have $L(n)=\lfloor\log_{{c}/{M}}
\frac{1}{n}\rfloor$ and the condition (\ref{eq:2.13}) is satisfied for any subsequence $(n_i)_i$  such that $
\log_{{c}/{M}}\frac{1}{n_i}-\lfloor\log_{c/M}\frac{1}{n_i}\rfloor$ 
converges (to $\log_{{c}/{M}}
\frac{1}{\kappa})$, and any $\kappa\in[1,\frac{M}c)$ can be obtained in this way.

\noindent
(b) As it will be seen in the proof, part (b) of the theorem holds under the weaker (than $r_{j+1}/r_j\to1$) assumption
\begin{equation}
\label{eq:2.17}
\lim_{j\to\infty}\frac{r_{j+1}}{h(j)}=0.
\end{equation}
(c) It is well known that Gaussian processes with covariances of the form (\ref{eq:2.14}) have the Markov property. This is in sharp contrast with the corresponding result in the Euclidean case [13,14].

\noindent
(d) The limit process obtained in part (b) has independent increments.

\noindent
(e) The function $R(t)$ accounts for the time scaling of the limit process only. From the form of the covariances it is seen that the most natural forms of $R(t)$ are
\begin{eqnarray}
\label{eq:2.18}
R(t)&=&\left\lfloor\log_a\frac{1}{t}\right\rfloor\,\,\hbox{\rm if}\,\, a<1,\\
\label{eq:2.19}
R(t)&=&\left\lfloor\log_M {t}\right\rfloor\,\,\hbox{\rm if}\,\, a=1.
\end{eqnarray}
We then have, for $a<1$,
\begin{equation}
\label{eq:2.20}
t^\theta\leq M^{R(t)}\leq Mt^\theta,
\end{equation}
where
\begin{equation}
\label{eq:2.21}
\theta=\frac{\log M}{\log\frac{1}{a}},
\end{equation}
and
\begin{equation}
\label{eq:2.22}
t\leq M^{R(t)}\leq Mt\,\,\,\,{\rm if}\,\, a=1.
\end{equation}
It is also seen that the exponents in (\ref{eq:2.15}) are ``close'' to $\frac{1}{t}$, namely
\begin{equation}
\label{eq:2.23}
\frac{1}{t}\leq a^{R(t)}<\frac{1}{at}.
\end{equation}
In (\ref{eq:2.20}), (\ref{eq:2.22}) and (\ref{eq:2.23}), the left-hand  inequalities are equalities for $t=a^{-n}, n\in \zet$.}
\end{remark}

Let us denote by $X^{(\kappa)}$ and $X$ the limit processes corresponding to $R$ given by (\ref{eq:2.18}) and (\ref{eq:2.19}), respectively. 
Observe that in spite of the fact that such $R$'s are not defined at $t=0$, the processes $X^{(\kappa)}$ and $X$ themselves are right continuous in $L^2$ at $0$ with $X^{(\kappa)}{(0)}=0$ and $X(0)=0$. This follows from (\ref{eq:2.20}), (\ref{eq:2.22}), and since the function $g_\kappa$ is bounded.

\begin{remark}
\label{r:2.10}
{\rm The process $X$ has the representation
$$X_t=\sqrt{1-e^{-2}}W_{M^n}\,\,\,\,{\rm for}\,\, M^n\leq t<M^{n+1},\ n\in\zet,$$
where $W$ is a standard Brownian motion.

Properties of the process $X^{(\kappa)}$ are summarized in the next 
proposition.}
\end{remark}

\begin{proposition}
\label{p:2.11}
(a) $X^{(\kappa)}$ is determined by a Gaussian sequence of random variables
$(\zeta_n)_{n\in\zetita}$,
$$X^{(\kappa)}_t=\zeta_n\,\,\,\,{ for}\,\, a^{-n}\leq t<a^{-(n+1)},$$
where the $\zeta_n$ have the representation
\begin{equation}
\label{eq:2.24}
\zeta_n=g_\kappa(a^{-n})
\sum^n_{i=-\infty}\nu_i
\sqrt{
\frac{M^i}{g_\kappa(a^{-i})}-
\frac{M^{i-1}}{g_\kappa(a^{-(i-1)})}},
\end{equation}
with $(\nu_i)_{i\in\zetita}$ i.i.d. standard normal.

\noindent
(b) $X^{(\kappa)}$ has a long range dependence
 property with exponent $1$ in the sense that
\begin{equation}
\label{eq:2.25}
0<\lim\sup_{\kern-.8cm\tau\to\infty}\tau|E(X^{(\kappa)}(t)-X^{(\kappa)}(s))(X^{(\kappa)}(t+\tau)-X^\kappa(s+\tau))|<\infty.
\end{equation}
(c) For $\theta>1$ (i.e., $1>a>\frac{1}{M}$), if we define
$$Y^{(\kappa)}_m(t)=a^{m(\theta-1)/2}X^{(\kappa)}(a^{-m}t),\,\,\,t\geq 0,$$
then
$$Y^{(\kappa)}_m\Rightarrow_fC(\kappa)
\sum_{n\in\zetita}Y(a^{-n})\UNO_{[a^{-n},a^{-(n+1)})}\,\,\,\,{ as}\,\, m\to\infty,$$
where
$$Y(t)=\frac{\sqrt{\theta+1}}{t}\int^t_0u^{\theta/2}dW_u,\,\,\, t\geq 0.$$
\end{proposition}

\begin{remark}
\label{r:2.12}
{\rm (a) Part (c) corresponds to Proposition 4.16 in [13]. Note that in that paper the limit process is a fractional Brownian motion, whereas in our case $Y$ 
is a diffusion process satisfying the equation
$$dY(t)=-\frac{1}{t}Y(t)dt +(\sqrt{\theta+1})t^{\theta/2-1}dW_t.$$
Its covariance has the particularly simple form
\begin{equation}
\label{eq:2.26}
EY(s)Y(t)=(s\wedge t)^\theta (s\vee t)^{-1}.
\end{equation}
(b) The condition $\theta>1$ in part (c) is required only to ensure that the process $Y$ is well defined at $0$. For $\theta\leq 1$ the assertion of part (c)
remains true on $[\varepsilon,\infty),\,\, \varepsilon >0$.

\noindent
(c) We think that part (a)  can be given a nicer form, analogous to that of the limit of $Y^{(\kappa)}_m$ in part (c), i.e., that the process $X^{(\kappa)}$ can be ``interpolated'' between the time points $a^{-n}$ by a Gaussian diffusion process. To derive such a representation it would suffice to prove that the function $g_\kappa$ defined by (\ref{eq:2.15}), with $a^{R(t)}$ replaced by 
$\frac{1}{t}$, is decreasing in $t$. We have not been able to prove it.}
\end{remark}

\begin{remark}
\label{r:2.13}
{\rm (a) Part (c) of Proposition 2.11 and Remark 2.12(b) show that the behavior of $X^{(\kappa)}(t)$ as $t\to\infty$ changes drastically depending on whether 
$\theta>1$ or $\theta<1$ (for $\theta>1, \lim_{t\to\infty}{\rm Var}X^{(\kappa)}(t)=\infty$, and for $\theta<1,\lim_{t\to\infty}{\rm Var}X^{(\kappa)}(t)=0$).
In fact, the parameter $\theta$ has a  deeper probabilistic interpretation which is related to the recurrence/transience properties of the underlying random walk. These properties are in a sense characterized by the degree $\gamma$ of a random walk  (see the Appendix). From (A.2) it follows that $\theta=\gamma+1$. In particular, this implies that for $\theta>1$ the random walk is transient, and for $\theta<1$ it is recurrent. In the critical case $\theta=1$ both situations can occur.
Moreover, by Theorem 2.2 it is not hard to check that if $\theta>1$, then 
$\lim_{L\to\infty}{\rm Var}N_n(L)=\infty$, and if $\theta<1$, then $\lim_{L\to\infty}{\rm Var} N_n(L)=0$, while, as it was seen in Example 2.3, if $\theta=1$ the number variance can be finite or infinite.
We recall, however, that in general the condition $\lim\sup_{L\to\infty}N_n(L)<\infty$ is not equivalent to recurrence of the random walk.

\noindent
(b) The equality $\theta=\gamma+1$ implies that for $c$-rw's the parameter $\theta$ corresponds to $\frac{1}{\alpha}$ for the symmetric $\alpha$-stable process in $\erre$ (see (\ref{eq:1.3}), and (A.2) in the Appendix). Hence for $c$-rw's the rate of the norming $\sqrt{M^{L(n)}}$ ($\sim n^{\theta/2}$) corresponds exactly to $t^{1/2\alpha}$, which is the norming in the fluctuation theorem in [13] (see (\ref{eq:1.5})).}
\end{remark}

Finally, we give
 the result for  the Poisson system.

\begin{theorem}
\label{t:2.14}
Under the assumptions of Proposition 2.5, for any $r_j$-rw, any functions 
$L:\zet_+\to\zet_+$ and $R:\erre_+\to\zet$, non-decreasing with 
$\lim_{n\to\infty}L(n)=\infty$ and $\lim_{t\to\infty}R(t)=\infty$, if $X_n$ is defined by (\ref{eq:2.11}), then $X_n\Rightarrow_f CX^{(R)}$, where $X^{(R)}$ is given in Theorem 2.8(b) and $C$ is a constant.
\end{theorem}

\begin{remark}
\label{r:2.15}
{\rm 
From the proof of this theorem it will be seen that the same result is true if the initial number of particles at site $u\in \Omega_M$ is $\eta^u$, where $(\eta^u)_{u\in\Omega_M}$ are i.i.d.\ random variables with $E(\eta^u)^{2+\delta}<\infty$ for some $\delta>0$, and $E\eta^u=\textrm{Var}\;\eta^u$.
}
\end{remark}

\section{Proofs}
\setcounter{section}{3}
\setcounter{equation}{0}

Let us denote
\begin{eqnarray*}
B_L(u)&=&\{x\in\Omega_M:|x-u|\leq L\},\,\,\,\, B_L=B_L(0),\\
S_L(u)&=&\{x\in\Omega_M:|x-u|= L\},\,\,\,\, S_L=S_L(0),\\
p^u_n(L)&=&P(u+\xi_n\in B_L).
\end{eqnarray*}
We write $|A|$ for the number of points in a bounded subset $A$ of $\Omega_M$ (no confusion arises with the hierarchical distance). Note that
$$|B_L(u)|=M^L,\quad |S_L(u)|=(M-1)M^{L-1}.$$
Observe that by ultrametricity of $\Omega_M$, if $u\in B_L$, then 
$u+\xi_n\in B_L$ if and only if $\xi_n\in B_L$. On the other hand, 
if $|u|=L+k, k=1,2,\ldots$, then $B_L\subset S_{L+k}(u)$, hence
$$p^u_n(L)=P(u+\xi_n\in B_L\cap S_{L+k}(u))=P(|\xi_n|=L+k)
\frac{|B_L|}{|S_{L+k}|}.$$
Thus we have
\begin{equation}
\label{eq:3.1}
p^u_n(L)=\left\{\begin{array}{lll}
P(|\xi_n|\leq L)&{\rm if}&u\in B_L,\\
P(|\xi_n|=L+k)\displaystyle\frac{1}{(M-1)M^{k-1}}&{\rm if}&|u|=L+k,\,\, k=1,2,\ldots .
\end{array}\right.
\end{equation}
{\bf Proof of Proposition 2.1}
By (\ref{eq:2.1a}) and (\ref{eq:3.1}),
$$EN_n(L)=\sum_{u\in B_L}p^u_n(L)+\sum_{u\notin B_L}p^u_n(L)=|B_L|.$$
This proves (a) (this formula is surely known, but we have not found a reference ).

In order to investigate Var$N_n(L)$ we study the tail of $\xi_n$. 
We have
\begin{eqnarray}
P(|\xi_n|>L)&\leq&\sum^n_{k=1}P(|\xi_1|\leq L,\ldots,|\xi_{k-1}|\leq L,|\xi_k|>L)\nonumber\\
\label{eq:3.2}
&=&\sum^n_{k=1}P(|\rho|\leq L)^{k-1}P(|\rho|>L),
\end{eqnarray}
by ultrametricity. Again by ultrametricity, the event that in the first $n$ steps there is just one farthest jump whose length is $j$ and it occurs at step $k$ has the form
$$A_{k,j}=\{|\xi_1|<j,\ldots,|\xi_{k-1}|<j,|\xi_k|=j,\,\,\xi_{k+1}\in B_{j-1}(\xi_k),\ldots,\xi_n\in B_{j-1}(\xi_k)\},$$
hence
\begin{eqnarray}
P(|\xi_n|>L)&\geq&\sum^n_{k=1}\sum^\infty_{j=L+1}P(A_{k,j})=\sum^n_{k=1}\sum^\infty_{j=L+1}P(|\rho|<j)^{n-1}P(|\rho|=j)\nonumber\\
\label{eq:3.3}
&\geq&n P(|\rho|\leq L)^{n-1}P(|\rho|>L).
\end{eqnarray}
Formulas (\ref{eq:3.2}) and (\ref{eq:3.3}) imply that
\begin{equation}
\label{eq:3.4}
P(|\xi_n|>L)\sim nP(|\rho|>L)\quad{\rm as}\quad L\to\infty.
\end{equation}
Using (\ref{eq:2.1a}) and (\ref{eq:3.1}) we have
$${\rm Var} N_n(L)=I(L)+I\!\!I(L),$$
where
\begin{eqnarray*}
I(L)&=&M^LP(|\xi_n|\leq L)P(|\xi_n|>L),\\
I\!\!I(L)&=&M^L\sum^\infty_{k=1}P(|\xi_n|=L+k)\left(1-\frac{1}{(M-1)M^{k-1}}P(|\xi_n|=L+k)\right).
\end{eqnarray*}
Hence (b) follows immediately from (\ref{eq:3.4}). 

The upper estimate in (c) is clear from the previous calculations. The lower bound is an easy consequence of
$${\rm Var} N_n(L)\geq I\!\!I(L)\geq M^LP(|\xi_n|>L)(1-P(|\xi_n|>L))$$
and (\ref{eq:3.2}), (\ref{eq:3.3}). \hfill$\Box$
\vglue.5cm
\noindent
{\bf Proof of Theorem 2.2} All the statements follow from 
(\ref{eq:2.2}) and the trivial formulas
$$M^LP(|\rho|>L)=M^Lr_{L+1}+M^LP(|\rho|>L+1),\,\,\,
M^LP(|\rho|>L)=\sum^\infty_{j=1}\frac{r_{j+L}M^{j+L}}{M^j}.$$
\hfill$\Box$
\vglue.5cm
\noindent
{\bf Proof of Proposition 2.5} If the number of particles at each site is Poisson with parameter $\lambda$, then
$${\rm Var} N_n(L)=\lambda \sum_{u\in\Omega_M} p^u_n(L)=\lambda EN_n(L)=\lambda M^L,$$
by (\ref{eq:2.1}), hence the assertion follows. \hfill$\Box$
\vglue.5cm
\noindent
{\bf Proof of Lemma 2.7} (a) By (\ref{eq:2.6}) and (\ref{eq:2.8}) it is easy to see that for any $0<\varepsilon <a\wedge (1-a)$ there exist positive constants 
$C_1, C_2$ such that for sufficiently large $n$,
\begin{equation}
\label{eq:3.5}
C_1(a-\varepsilon)^n\leq h(n)\leq C_2(a+\varepsilon)^n.
\end{equation}
Then, for large $n$, 
$$h\left(\left\lfloor \log_{a+\varepsilon}\frac{1}{n}+\log_{a+\varepsilon}\frac{1}{C_2}\right\rfloor +1\right)\leq C_2(a+\varepsilon)
^{\log_{(a+\varepsilon)}{1/nC_2}}=\frac{1}{n},$$
which implies that
$$L(n)\leq \log_{a+\varepsilon}\frac{1}{n}+\log_{a+\varepsilon}\frac{1}{C_2}+1,$$
hence
$$\lim\sup_{\kern-.8cm n\to\infty}\frac{L(n)}{\log_a\frac{1}{n}}\leq \frac{\log a}{\log (a+\varepsilon )}.$$
Analogously, using the left-hand side of (\ref{eq:3.5}) we obtain
$$\lim\inf_{\kern-.8cm n\to\infty}\frac{L(n)}{\log_a\frac{1}{n}}\geq \frac{\log a}{\log (a-\varepsilon )}.$$
Hence (\ref{eq:2.9}) follows

\noindent
(b) The lower bound in (\ref{eq:2.10}) is obvious by (\ref{eq:2.7}). The assumption
(\ref{eq:2.6}) clearly implies
\begin{equation}
\label{eq:3.6}
\lim_{n\to\infty}\frac{r_{n+k}}{r_n}=a^k,\qquad k=1,2,\ldots,
\end{equation}
hence for $a<1$ we have
$$\left(\sum^\infty_{j=0}(a+\varepsilon)^j\right)^{-1}\leq
\lim\inf_{\kern-.8cm n\to\infty}\frac{r_{n+1}}{h(n)}\leq 
\lim\sup_{\kern-.8cm n\to\infty}\frac{r_{n+1}}{h(n)}
\leq\left(\sum^{k}_{j=0}a^j\right)^{-1}$$
for any $0<\varepsilon <1-a$ and $k=1,2,\ldots$. The right-hand side inequality also holds for $a=1$. Therefore, for all $0<a\leq 1$,
\begin{equation}
\label{eq:3.7}
\lim_{n\to\infty}\frac{r_{n+1}}{h(n)}=1-a.
\end{equation}

By (\ref{eq:2.7}) and (\ref{eq:2.8}),
$$nh(L(n))<1+nr_{L(n)+1},$$
so (\ref{eq:2.10}) follows easily from (\ref{eq:3.7}). \hfill$\Box$
\vglue.5cm
\noindent
{\bf Proof of Theorem 2.8} (a) We will prove
\begin{equation}
\label{eq:3.8}
\lim_{i\to\infty}{\rm Cov}(X_{n_i}(s), X_{n_i}(t))=M^{R(s\wedge t)}g_\kappa
(s\vee t)
\end{equation}
and 
\begin{equation}
\label{eq:3.9}
\lim_{n\to\infty}\frac{1}{(M^{L{(n)}})^{1+{\delta}/{2}}}\sum_{u\in\Omega_M}
E\left|\UNO_{B_{(L(n)+R(t))^+}}(u+\xi^u_n)-P(u+\xi^u_n\leq(L(n)+R(t))^+)\right|^{2+\delta}=0,\quad \delta>0.
\end{equation}
Then by (\ref{eq:2.1a}), (\ref{eq:2.11}) and independence of random walks the result follows from the central limit theorem in the Lyapunov version.

Denote
\begin{equation}
\label{eq:3.10}
L(i,t)=(L(n_i)+R(t))^+.
\end{equation}
Fix $s\leq t$. By (\ref{eq:2.1a}) and independence of random walks,
\begin{eqnarray}
{\rm Cov}(X_{n_i}(s),X_{n_i}(t))&=&M^{-L(n_i)}
\sum_{u\in\Omega_M}p^u_{n_i}(L(i,s))\left(1-p^u_{n_i}(L(i,t))\right)\nonumber\\
\label{eq:3.11}
&=&M^{-L(n_i)}\left(I+I\!\!I\right),
\end{eqnarray}
where $I=\sum_{u\in B(L(i,t))}\ldots$ and $I\!\!I=\sum_{u\notin B(L(i,t))}\ldots$.
We have
\begin{eqnarray}
I&=&\sum_{u\in B(L(i,s))}\ldots +\sum_{u\in B(L(i,t))\backslash 
B(L(i,s))}\ldots\nonumber\\
&=&M^{L(i,s)}P(|\xi_{n_i}|\leq L(i,s))P(|\xi_{n_i}|>L(i,t))\nonumber\\
&&+\sum^{L(i,t)}_{j=L(i,s)+1}\sum_{u\in S_j}P(|\xi_{n_i}|=j)\frac{|B_{L(i,s)}|}{|S_j|}\left(1-P(|\xi_{n_i}|\leq L(i,t))\right)\nonumber\\
\label{eq:3.12}
&=&M^{L(i,s)}P(|\xi_{n_i}|\leq L(i,t))P(|\xi_{n_i}|>L(i,t)),
\end{eqnarray}
by (\ref{eq:3.1}). Similarly
\begin{eqnarray}
I\!\!I&=&\sum^\infty_{j=1}\sum_{u\in S_{L(i,t)+j}}\ldots\nonumber\\
\label{eq:3.13}
&=&M^{L(i,s)}\sum^\infty_{j=1}P(|\xi_{n_i}|=L(i,t)+j)\left(1-P(|\xi_{n_i}|=L(i,t)+j)\frac{1}{(M-1)M^{j-1}}\right).
\end{eqnarray}
By (\ref{eq:3.10})-(\ref{eq:3.13}), for sufficiently large $i$, we obtain
\begin{eqnarray}
{\rm Cov}(X_{n_i}(s),X_{n_i}(t))&=&M^{R(s)}
\biggl[ P(|\xi_{n_i}|>L(i,t))\biggr.
(1+P(|\xi_{n_i}|\leq L(i,t)))\nonumber\\
\label{eq:3.14}
&&-\biggl.\frac{1}{M-1}\sum^\infty_{j=1}\frac{P^2(|\xi_{n_i}|=L(i,t)+j)}
{M^{j-1}}\biggr].
\end{eqnarray}
It is known ([20], see also [7]) that
$$P(\xi_n=u)=-\frac{f^n_k}{M^k}+(M-1)
\sum^\infty_{j=k+1}\frac{f^n_j}{M^j},\quad{\rm if}\quad |u|=k>0,$$
where
\begin{equation}
\label{eq:3.15}
f_k=1-h(k-1)-\frac{r_k}{M-1}.
\end{equation}
Hence
\begin{eqnarray}
P(|\xi_n|=k)&=&\frac{M-1}{M}\left(-f^n_k+(M-1)\sum^\infty_{j=1}\frac{f^n_{j+k}}{M^j}\right)\nonumber \\
\label{eq:3.16}
&=&\frac{M-1}{M}\sum^\infty_{j=0}\frac{f^n_{j+k+1}-f^n_{j+k}}{M^j}.
\end{eqnarray}
This implies that
\begin{equation}
\label{eq:3.17}
P(|\xi_n|>L)=\frac{M-1}{M}\sum^\infty_{j=0}\frac{1}{M^j}(1-f^n_{L+j+1}).
\end{equation}
By (\ref{eq:2.13}), (\ref{eq:3.6}), (\ref{eq:3.7}) and (\ref{eq:3.10}) we have
\begin{equation}
\label{eq:3.18}
\lim_{i\to\infty}n_i\left(h(L(i,t)+j-1)+\frac{r_{L(i,t)+j}}{M-1}\right)=\kappa b a^{R(t)+j-1},
\end{equation}
where $b$ is defined by (\ref{eq:2.12}). Hence
\begin{equation}
\label{eq:3.19}
\lim_{i\to\infty}f^{n_i}_{L(i,t)+j}=e^{-\kappa b a^{R(t)+j-1}},
\end{equation}
by (\ref{eq:3.15}).

Combining (\ref{eq:3.16})-(\ref{eq:3.19}) and (\ref{eq:3.14}) (it is clear that we can pass to the limits under the sums) we arrive at
\begin{eqnarray}
\lefteqn{\lim_{i\to\infty}{\rm Cov}(X_{n_i}(s),X_{n_i}(t))}\nonumber\\
&=&M^{R(s)}\left[\left(\frac{M-1}{M}\sum^\infty_{j=0}
\frac{1-e^{-\kappa b a^{R(t)}a^j}}{M^j}\right)\left(1+
\frac{M-1}{M}\sum^\infty_{j=0}
\frac{e^{-\kappa b a^{R(t)}a^j}}{M^j}\right)\right.\nonumber\\
&&-\frac{(M-1)^3}{M^4}\sum^\infty_{j=0}
\frac{1}{M^j}\left(\sum^\infty_{k=0}
\left. 
\frac{e^{-\kappa b a^{R(t)}a^{j+1+k}}-e^{-\kappa b a^{R(t)}a^j}}
{M^k}\right)^2\right]\nonumber\\
\label{eq:3.20}
&=&M^{R(s)}g_\kappa(t)
\end{eqnarray}
(see (\ref{eq:2.15})). This proves (\ref{eq:3.8}).

It is easy to see that the expression under the limit on the left hand side of 
(\ref{eq:3.9}) is estimated from above by
$$\frac{C}{M^{L(n)(1+\delta/2)}}\sum_{u\in\Omega_M}p^u_n((L(n)+R(t))^+)=\frac{C}{M^{L(n)(1+\delta/2)}}M^{(L(n)+R(t))^+},$$
by (\ref{eq:2.1}), hence (\ref{eq:2.9}) follows.

The proof of part (a) is complete.

\noindent
(b) The fact that 
\begin{equation}
\label{eq:3.21}
\lim_{n\to\infty}n(h(L(n))=1
\end{equation}
follows immediately from (\ref{eq:2.10}). This and (\ref{eq:3.7}) imply that 
\begin{equation}
\label{eq:3.22}
\lim_{n\to\infty}n\left(h((L(n)+R(t))^++j-1)+\frac{r_{(L(n)+R(t))^++j}}{M-1}\right)=1
\end{equation}
(cf. (\ref{eq:3.18})), hence
$$\lim_{n\to\infty}f^n_{(L(n)+R(t))^++j}=e^{-1}$$
for $j=0,1,2,\ldots$. Therefore, a counterpart of (\ref{eq:3.20}) is 
\begin{equation}
\label{eq:3.23}
\lim_{n\to\infty}{\rm Cov}(X_n(s),X_n(t))=(1-e^{-2})M^{R(s\wedge t)}.
\end{equation}
Now, the $\Rightarrow_f$ convergence (the Lyapunov condition (\ref{eq:3.9})) is obtained in the same way as before.

Observe that (\ref{eq:3.21}), (\ref{eq:3.22}), and hence (\ref{eq:3.23}) as well, follow from (\ref{eq:3.7}) in this case. (Recall that  (\ref{eq:3.7}) with $a=1$ implies (\ref{eq:2.10}), see the end of the proof of Lemma \ref{l:2.7}(b).)
Therefore, as stated in Remark 2.9 (b), assumption (\ref{eq:2.6}) can be replaced by (\ref{eq:2.17}). \hfill$\Box$
\vglue.5cm
\noindent
{\bf Proof of Proposition 2.11} (a) From positive-definiteness of the covariance function (\ref{eq:2.14}) it is easy to see that 
$M^i/g_\kappa {(a^{-i})}$ is increasing in $i$, moreover, by (\ref{eq:2.15}), 
$\lim_{i\to-\infty}M^i/g_\kappa{(a^{-i}})=0$. Hence the result is obtained by a direct computation.

\noindent
(b) By (\ref{eq:2.14}), for $s<t<s+\tau<t+\tau$ we have 
\begin{eqnarray*}
\lefteqn{E(X^{(\kappa)}(t)-
X^{(\kappa)}(s))(X^{(\kappa)}(t+\tau)-X^{(\kappa)}(s+\tau))}\\
&=&(M^{R(t)}-M^{R(s)})(g_\kappa(t+\tau)-g_\kappa(s+\tau)),
\end{eqnarray*}
hence it suffices to investigate the second factor. Obviously, for large $\tau$, the time points $s+\tau, t+\tau$ belong either to the same interval of the form $[a^{-k}, a^{-(k+1)})$, or to two neighboring such intervals. Since, by 
(\ref{eq:2.18}) and (\ref{eq:2.15}), $g_\kappa$ is constant on such intervals, it is enough to consider a sequence of the form
$$\tau_m=a^{-m}-d_m,\quad s<d_m\leq t,\quad m=1,2,\ldots.$$
Then, for large $m$ we have
$$\tau_m(g_\kappa(\tau_m+t)-g_\kappa(\tau_m +s))=\tau_ma^{m}a^{-m}(\tilde{g}_\kappa(a^m)-\tilde{g}_\kappa(a^{m-1})),$$
where $\tilde{g}_\kappa(r)$ is obtained from (\ref{eq:2.15}) by putting $r$ instead of $a^{R(t)}$. The mean value theorem implies that $\frac{1}{r}(\tilde{g}_\kappa(r)-\tilde{g}_\kappa(\frac{r}{a}))$ has a finite positive limit as $r\to 0$. This proves (\ref{eq:2.25}).

\noindent
(c) Since the processes $Y^{(\kappa)}_m$ are centered Gaussian, it suffices to prove convergence of covariances. By (\ref{eq:2.18}) and (\ref{eq:2.21}),
\begin{eqnarray}
\label{eq:3.24}
R(a^{-m}t)&=&m+R(t),\quad m=1,2,\ldots,\\
M^m&=&a^{-m\theta}.\nonumber
\end{eqnarray}
This, (\ref{eq:2.14}) and (\ref{eq:2.18}) imply that for $s\leq t$,
$${\rm Cov}(Y^{(\kappa)}_m(s), Y^{(\kappa)}_m(t))=M^{R(s)}
a^{-m}g_\kappa(a^{-m}t).$$
By (\ref{eq:3.24}) we also have
$$\lim_{m\to\infty}\frac{1-e^{-\kappa b a^{R(a^{-m}t)}a^j}}{a^m}=\kappa b a^{R(t)}a^j.$$
Hence, using the form of $g_\kappa$ given by (\ref{eq:3.20}), and by 
(\ref{eq:2.12}), it is easy to see that 
$$\lim_{m\to\infty}{\rm Cov}(Y^{(\kappa)}_m(s), Y^{(\kappa)}_m(t))=2\kappa M^{R(s)}a^{R(t)}, \quad s\leq t.$$
Hence (c) follows (cf (\ref{eq:2.26})). \hfill$\Box$
\vglue.5cm
\noindent
{\bf Proof of Theorem 2.14}
Let $\eta^u$ be the number of particles at site $u$ at time $0$. The 
{r.v.'s} $(\eta^u)_{u\in\Omega_M}$ are i.i.d.
Let $(\xi^{u,k})_{u\in\Omega_M,k\in\enita}$ be independent copies of 
$r_j$-${\rm rw}$.
Then
$$N_n(L)=\sum_{u\in\Omega_M}\sum^{\eta^u}_{k=1}\UNO_{\{u+\xi^{u,k}_n\in 
B_L\}},$$
and for $s\leq t$, denoting $L(n,t)=(L(n)+R(t))^+$,
\begin{eqnarray*}
{\rm Cov}(X_n(s),X_n(t))&=&\frac{1}{M^{L(n)}}
\sum_{u\in\Omega_M}\left[E\eta^u P(u+\xi_n\in B_{L(n,s)}\right.)\\
&&+\left.({\rm Var}\,\eta^u-E\eta^u)P(u+\xi_n\in B_{L(n,s)})P(u+\xi_n\in B_{L(n,t)})\right]\\
&=&\frac{E\eta^0}{M^{L(n)}}M^{L(n,s)},
\end{eqnarray*}
by (\ref{eq:2.1}). Hence
$$\lim_{n\to\infty}{\rm Cov}(X_n(s), X_n(t)=E\eta^0 M^{R(s)}.$$
Next, using the inequality $(a_1+\ldots+a_m)^{2+\delta}\leq m^{1+\delta}(a^{2+\delta}_1+\ldots+a^{2+\delta}_m), a_i\geq 0$, we easily obtain
$$\lim_{n\to\infty}\frac{1}{M^{(1+\delta/2)L(n)}}\sum_{u\in\Omega_M}E
\left(\sum^{\eta^u}_{k=1}\UNO_{B_{L(n,t)}}(u+\xi^{u,k}_n)\right)^{2+
\delta}=0.$$
Convergence of finite-dimensional distributions now follows from the central limit theorem (Lyapunov criterion). \hfill$\Box$

\section*{Appendix}
\setcounter{equation}{0}
\def\theequation{{\rm A}.\arabic{equation}}

\vglue.5cm
The development and applications of hierarchical random walks are outlined in \cite{DGW3}, and their recurrence/transience  and some other properties are studied in \cite{DGW2} (and references therein).

In order to deal with the standard $\alpha$-stable process on $\erre^d$ and discrete time hierarchical random walks in a unified way, the continuous time version of a walk with unit rate holding time was taken in \cite{DGW2}. Its transition probability from $0$ to $y$ in time $t>0$ is given by
$$p_t(0,y)=e^{-t}\sum^\infty_{n=0}\frac{t^n}{n!}p^{(n)}(0,y),$$
where $p^{(n)}(0,y)$ is the transition probability from $0$ to $y$ in $n$ steps. Generally, for a L\'evy process on a Polish space $S$ (with additive group structure in the cases considered in \cite{DGW2}) with semigroup ${\cal T}_t$, the {\it degree} of the process is defined by 
$$\gamma=\sup\{\zeta >-1:G^{\zeta+1}\varphi<\infty\quad\hbox{\rm for all}\quad \varphi\in{\cal B}^+_b(S)\},$$
where ${\cal B}^+_b(S)$ is the space of bounded non-negative measurable functions with bounded support on $S$, and
$$G^\zeta\varphi=\frac{1}{\Gamma(\zeta)}\int^\infty_0t^{\zeta-1}{\cal T}_t\varphi dt,\quad\zeta>0,$$
is the fractional power of the potential (or Green) operator of the process. For $\zeta=\gamma$ it may happen that $G^{\gamma+1}\varphi<\infty$ or $G^{\gamma+1}\varphi=\infty,\varphi\neq 0$.   Recurrence corresponds to $\gamma\in (-1,0)$ and transience corresponds to $\gamma\in (0,\infty]$. The value $\gamma=0$ is special since both cases $G\varphi<\infty$ and 
$G\varphi=\infty,\,\,\varphi\neq 0$, can happen. $\gamma=\infty$ occurs, for example, for a simple asymmetric random walk on $\zet$. For $\gamma \geq 0$, $\gamma$ is also given by
$$
\gamma =\sup \{ \zeta \geq 0:\,\,EL^\zeta_{B_R}<\infty\quad \hbox{\rm for all}\quad R>0\},
$$
where $L_{B_R}$ is the last exit time of the process (starting at $0$) from an open ball $B_R$ of radius $R$ centered at $0$. For a discrete time random walk on a discrete space, instead of $G^\zeta\varphi$ it is natural to consider
\begin{equation}
\label{A.1}
\frac{1}{\Gamma (\zeta)} \displaystyle\sum^\infty_{n=0} n^{\zeta-1} p^{(n)}(0,0),
\end{equation}
\noindent
but since
$$
\int^\infty_0 t^\zeta p_t (0,0) dy =\sum^\infty_{n=0} \frac{\Gamma (n+\zeta +1)}{n!} p^{(n)} (0,0),
$$
and
$$
\frac{\Gamma (n+\zeta +1)}{n!n^{\zeta}}\rightarrow 1 \quad {\rm as}\quad n\rightarrow \infty,
$$
the degree of the random walk is the same in discrete and continuous time.

For the standard $\alpha$-stable process on $\erre^d, \gamma$ is given by
(\ref{eq:1.3}),
and the range of possible values of $\gamma$ for those processes is  restricted to the interval $[-\frac{1}{2},\infty)$. For the $c$-rw  defined by (\ref{eq:1.1}),(\ref{eq:1.2}), $\gamma$ is given by (\ref{eq:1.4}), 
and the range of possible values of $\gamma$ 
is $(-1,\infty)$. Hence this class of hierarchical random walks is richer than the class of standard $\alpha$-stable processes and corresponding symmetric random walks on $\zet^d$. The $c$-rw's mimic the behavior of $\alpha$-stable processes in the sense that they have the same recurrence/transience properties for equal values of their degrees; in particular, for $\gamma=0, \,\,G\varphi=\infty$ holds for both of them. The correspondence of degrees allows to choose $c$ in order to study ``caricatures'' of $\alpha$-stable processes by means of 
$c$-rw's, including non-integer values of the dimension $d$ (which is one of the reasons for using hierarchical random walks in statistical physics, see 
\cite{DGW3,DGW4} and references therein).
But there are also differences, for example, the distance from 0 of the 
$c$-rw has a different behavior from the $\alpha$-stable Bessel process (see \cite{DGW2}, Remark 3.5.7). The correspondence between the degrees of the $\alpha$-stable process and the $c$-rw also serves to compare our results  with those of the $\alpha$-stable process  and exhibit analogies and differences.

The assumption  
 $r_{j+1}/r_j\to a$ as $j\to\infty, 0<a\leq 1$, implies that the degree of the
$r_j$-rw is given by
\begin{equation}
\label{A.2}
\gamma=\frac{\log M}{\log \frac{1}{a}}-1\quad{\rm for}\quad a<1,\,\,
\gamma=\infty\quad{\rm for}\quad a=1.
\end{equation}
This is obtained from \cite{DGW2} (Proposition 3.2.7). Hence
$$\gamma\left\{\begin{array}{ll}
<0,&\\
=0,&\hbox{\rm if}\\
>0,\end{array}\right.a\left\{\begin{array}{l}
<\frac{1}{M},\,\, {\rm recurrent}\\
=\frac{1}{M},\\
>\frac{1}{M},\,\, {\rm transient}.\end{array}\right.
$$
The parameter $\theta=\gamma+1$ (see (\ref{eq:2.21})) therefore has a meaning in terms of 
recurrence/transience.

Examples:

\begin{enumerate}
\item {\it c-rw}: $0<c<M$,
$$r_j=\left(1-\frac{c}{M}\right)\left(\frac{c}{M}\right)^{j-1},\quad j=1,2,\ldots,\quad a=\frac{c}{M}, \quad\gamma=\frac{\log c}{\log(M/c)},$$
recurrent for $c\leq 1$, transient for $c>1$.

\item $j^\beta$-{\it rw}:$\beta\geq 0$, 
$$r_j=D\frac{j^\beta}{M^j},\quad j=1,2,\ldots,$$
where $D$ is a normalizing constant,
$$a=\frac{1}{M},\quad \gamma=0,$$
recurrent for $\beta\leq 1$, transient for $\beta > 1$ (follows from the recurrence criterion in \cite{FP}, see also Example 3.2.6 in \cite{DGW2}, which is slightly different) .

\item 
$$r_j=Dj^{-\beta},\quad j=1,2,\ldots,$$
where $\beta >1$ and $D$ is a normalizing constant,
$$a=1,\quad \gamma=\infty,\quad{\rm transient}.$$
\end{enumerate}
{\bf Remark}. In \cite{DGW2} $r_j$ is written in the form
$$
r_j =D \frac{c_j}{N^{j/\mu}},\quad j=1,2,\ldots,
$$
where $\mu$ is a positive constant, $(c_j)$ is a sequence of positive numbers, and $D$ is a normalizing constant. The parameter $\mu$ is useful for some applications (see also \cite{DGW3,DGW4}).

\end{document}